\documentclass[a4paper,11pt]{article}
\usepackage[T1]{fontenc}
\usepackage{lmodern}
\usepackage[margin=27mm]{geometry}
\usepackage{amsmath,amssymb,amsthm,mathtools,booktabs,microtype,enumitem}
\usepackage[colorlinks=true,linkcolor=blue,citecolor=blue,urlcolor=blue]{hyperref}
\setlist[enumerate]{itemsep=3pt,topsep=5pt}
\newtheorem{theorem}{Theorem}[section]
\newtheorem{lemma}[theorem]{Lemma}
\newtheorem{proposition}[theorem]{Proposition}
\newtheorem{corollary}[theorem]{Corollary}
\theoremstyle{definition}
\newtheorem{definition}[theorem]{Definition}
\newtheorem{example}[theorem]{Example}
\newtheorem{remark}[theorem]{Remark}
\DeclareMathOperator{\ord}{ord}
\DeclareMathOperator{\rad}{rad}

\DeclareMathOperator{\wt}{wt}

\DeclareMathOperator{\Tr}{Tr}
\newcommand{\F}{\mathbb F}
\newcommand{\T}{\mathcal T}
\newcommand{\E}{\mathcal E}
\newcommand{\R}{\mathcal R}
\newcommand{\doi}[1]{\href{https://doi.org/#1}{doi: #1}}
\hypersetup{pdftitle={Cyclotomic polynomials whose irreducible factors are trinomials over finite fields},pdfauthor={},pdfkeywords={finite fields, cyclotomic polynomials, trinomials}}
\title{Cyclotomic polynomials with trinomial irreducible factors over finite fields}
\author{Li Zhu$^{1}$ and Hongfeng Wu$^{2}$\footnote{Corresponding author.}
	\setcounter{footnote}{-1}
	\footnote{E-Mail addresses:
		lizhumath@pku.edu.cn (L. Zhu), whfmath@gmail.com (H. Wu)}
	\\[0.5ex]
	\small $^{1}$School of Mathematical Sciences, Guizhou Normal University, Guiyang, China
	\\
	\small $^{2}$College of Science, North China University of Technology, Beijing, China}
\date{}
\begin{document}
\maketitle
\begin{abstract}
We classify the pairs $(q,n)$, where $q$ is a prime power, $n>1$, and $\gcd(n,q)=1$, for which every irreducible factor of $\Phi_n(x)$ over $\F_q$ has exactly three nonzero terms. The positive support exponents of these factors have a common greatest common divisor, which gives a unique reduction by a power substitution. The reduced degree is two in odd characteristic and belongs to $\{2,3,4\}$ in characteristic two; the cubic and quartic cases arise by scalar substitution from the layers of orders $7$ and $15$. This yields necessary and sufficient arithmetic conditions on $(q,n)$ and explicit factorizations. Among the classified families, precisely those arising from quartic bases have a divisor layer containing an irreducible factor with more than three terms.
\end{abstract}

\section{Introduction and main results}
If $\gcd(n,q)=1$, every irreducible factor of $\Phi_n(x)$ over $\F_q$ has degree $\ord_n(q)$, and all of its roots have multiplicative order $n$. These facts do not determine which coefficients vanish. We classify the pairs $(q,n)$ for which every factor in this cyclotomic layer has exactly three nonzero terms.

For a polynomial $f$, write $\wt(f)$ for the number of its nonzero coefficients. Throughout the paper, a \emph{trinomial} means a polynomial with exactly three nonzero terms. Thus a monic trinomial with nonzero constant coefficient has the form
\[
 x^d+A x^k+B,\qquad 0<k<d,\quad A,B\in\F_q^*.
\]
We write $\T_q(n)$ if every monic irreducible factor of $\Phi_n(x)$ over $\F_q$ is a trinomial. Unless otherwise stated, $n>1$ and $\gcd(n,q)=1$. We use $\rad(n)=\prod_{\ell\mid n}\ell$, with $\rad(1)=1$, and put $\ord_1(q)=1$ when this convention is needed.

Explicit sparse factorizations of cyclotomic polynomials and of $x^n-1$ were obtained in \cite{WW,BGO}; general formulas for cyclotomic polynomials and polynomial compositions are given in \cite{Graner}. Oliveira and Reis \cite{OR} studied the requirement that every irreducible factor of $x^n-1$ be a binomial or a trinomial over binary fields. Cheng \cite{ChengI,ChengII,ChengCorr} treated its classification over arbitrary finite fields. Citations to \cite{ChengII} below refer to the corrected author version listed in the bibliography.

Our condition concerns a single cyclotomic layer in the decomposition $x^n-1=\prod_{e\mid n}\Phi_e(x)$, rather than all divisor layers. For example, over $\F_2$,
\[
 \Phi_{15}(x)=(x^4+x+1)(x^4+x^3+1),
\]
whereas $\Phi_5(x)=x^4+x^3+x^2+x+1$ is irreducible. Corollary~\ref{cor:global} identifies all classified families for which this distinction occurs.

\begin{theorem}[Root-order classification]\label{thm:main}
Let $q=p^r$, where $p$ is prime, and let $n>1$ with $\gcd(n,q)=1$.
\begin{enumerate}[label=\textup{(\roman*)}]
\item If $p$ is odd, put $c=\gcd(n,q^2-1)$. Then
\begin{equation}\label{eq:odd-main}
 \T_q(n)\quad\Longleftrightarrow\quad
 \rad(n)\mid q^2-1\quad\text{and}\quad c\nmid 2(q-1).
\end{equation}
\item If $p=2$, then $\T_q(n)$ holds if and only if exactly one of the following alternatives holds:
\begin{align}
\textup{(Q)}\quad&\rad(n)\mid q^2-1,\qquad \rad(n)\nmid q-1;\label{eq:class2}\\
\textup{(C)}\quad&3\nmid r,\quad n=7^aM,\quad a\ge1,\quad\gcd(M,7)=1,\quad\rad(M)\mid q-1,\notag\\[-1mm]
 &3\nmid M/\gcd(M,q-1),\qquad a=1\ \text{or}\ 7\nmid r;\label{eq:class3}\\
\textup{(F)}\quad&r\text{ is odd},\quad n=3^a5^bM,\quad a,b\ge1,\quad\gcd(M,15)=1,\notag\\[-1mm]
 &\rad(M)\mid q-1,\qquad
 (a=1\ \text{or}\ 3\nmid r),\quad(b=1\ \text{or}\ 5\nmid r).
 \label{eq:class4}
\end{align}
\end{enumerate}
\end{theorem}
The labels \textup{(Q)}, \textup{(C)}, and \textup{(F)} refer to the degrees two, three, and four before a power substitution.

\begin{definition}\label{def:base}
A positive integer $e>1$ is a \emph{base order} over $\F_q$ if $\T_q(e)$ holds and some irreducible factor $x^D+A x^K+B$ of $\Phi_e$ satisfies $\gcd(D,K)=1$. Let $\E_q$ be the set of base orders and let $\E_q^{(j)}=\{e\in\E_q:\ord_e(q)=j\}$.
\end{definition}

\begin{theorem}[Base orders]\label{thm:bases}
If $q$ is odd, then
\begin{equation}\label{eq:odd-bases}
 \E_q=\E_q^{(2)}=\{e:e\mid q^2-1,\ e\nmid2(q-1)\}.
\end{equation}
If $q=2^r$, then $\E_q$ is the disjoint union of
\begin{align}
 \E_q^{(2)}&=\{e:e\mid q^2-1,\ e\nmid q-1\},\label{eq:base2}\\
 \E_q^{(3)}&=\begin{cases}\{7s:s\mid q-1\},&3\nmid r,\\ \varnothing,&3\mid r,\end{cases}\label{eq:base3}\\
 \E_q^{(4)}&=\begin{cases}\{15s:s\mid q-1\},&r\text{ odd},\\ \varnothing,&r\text{ even}.\end{cases}\label{eq:base4}
\end{align}
\end{theorem}

Proposition~\ref{prop:descent} reduces every admissible layer uniquely to a base order by a power substitution. The bounds on base degrees follow from power maps on roots, quadratic and cubic coefficient norms, and Vishne's extension \cite{Vishne} of Swan's parity method \cite{Swan}. In characteristic two, two polynomial identities exclude the remaining degree-twelve case. Section~\ref{sec:assembly} derives the lifting conditions and explicit factors, and Section~\ref{sec:consequences} gives examples and quartic factor counts.

\section{Frobenius orbits and power substitutions}\label{sec:reduction}
We use the standard Frobenius description of finite-field extensions \cite[Chapter~2]{LN}. If $\alpha$ has multiplicative order $n$, its conjugates over $\F_q$ are $\alpha,\alpha^q,\ldots,\alpha^{q^{d-1}}$, where $d=\ord_n(q)$. Thus $\Phi_n$ is squarefree and has $\varphi(n)/d$ irreducible factors. Over $\F_{q^t}$, a degree-$d$ irreducible polynomial splits into $\gcd(d,t)$ factors of degree $d/\gcd(d,t)$, by replacing Frobenius with its $t$th power. We write $\mu_m$ for the group of $m$th roots of unity in an algebraic closure of $\F_q$, and $\ord(\alpha)$ for the multiplicative order of a nonzero element. The monic reciprocal $f^*(x)=x^{\deg f}f(x^{-1})/f(0)$ preserves irreducibility, weight, and root order.

\subsection{A common support divisor}
For a nonconstant polynomial $f$ with $f(0)\ne0$, define
\[
 h(f)=\gcd\{i>0:[x^i]f\ne0\}.
\]
\begin{lemma}[Common support divisor]\label{lem:common-support}
The number $h(f)$ is the same for every monic irreducible factor $f$ of $\Phi_n(x)$ over $\F_q$.
\end{lemma}
\begin{proof}
Choose a primitive $n$th root $\zeta$ and let $H=\langle q\rangle$ in $(\mathbb Z/n\mathbb Z)^*$. The exponents of the roots of each irreducible factor form a coset $aH$, where $a$ is a unit modulo $n$. Their additive stabilizers are
\[
 S_a=\{b\in\mathbb Z/n\mathbb Z:aH+b=aH\}.
\]
Each $S_a$ is an additive subgroup and $S_a=aS_1$. Every subgroup of a cyclic group is preserved by multiplication by a unit, so $S_a=S_1$ for all $a$.

Let $f$ be a monic factor of degree $d$. A scalar $\eta$ preserves its root set if and only if
\[
 f(\eta x)=\eta^d f(x).
\]
Comparison of the nonzero constant terms first gives $\eta^d=1$; comparison of all other nonzero coefficients then gives $\eta^i=1$ for every positive support exponent $i$. Hence the scalar stabilizer is $\mu_{h(f)}$. Since $f$ is separable, the characteristic does not divide $h(f)$, so this stabilizer has exactly $h(f)$ elements. Moreover, every such scalar is a ratio of two roots of $f$ and therefore belongs to $\mu_n$. Under $\zeta^b\leftrightarrow b$, it corresponds exactly to $S_a$. The common size of the $S_a$ proves the assertion.
\end{proof}

\begin{lemma}\label{lem:cyclosub}
Let $e,h$ be positive integers. If $\rad(h)\mid e$, then $\Phi_{eh}(x)=\Phi_e(x^h)$.
\end{lemma}
\begin{proof}
Apply the identity $\Phi_{N\ell}(x)=\Phi_N(x^\ell)$ for $\ell\mid N$ successively to the prime factors of $h$, all of which divide $e$; see \cite{LN,WW}.
\end{proof}

\begin{proposition}[Canonical reduction]\label{prop:descent}
Suppose $\T_q(n)$ holds. For any irreducible factor
\[
 f=x^d+A x^k+B,\qquad h=\gcd(d,k),\quad D=d/h,\quad K=k/h,
\]
put $F=x^D+A x^K+B$. Then $h$ is independent of $f$ and
\begin{equation}\label{eq:descent}
 h\mid n,\quad e=n/h,\quad\rad(h)\mid e,\quad
 D=\ord_e(q),\quad e\in\E_q,\quad\ord_n(q)=h\ord_e(q).
\end{equation}
The irreducible factors of $\Phi_n$ are exactly the polynomials $G(x^h)$, where $G$ runs through the irreducible factors of $\Phi_e$. Every such $G$ has support divisor one. In particular, the pair $(e,h)$ is unique.
\end{proposition}
\begin{proof}
Lemma~\ref{lem:common-support} proves independence of $h$. Separability gives $p\nmid h$. If $\beta$ is a root of $f$ and $\xi$ is a primitive $h$th root, then $\xi\beta$ is also a root because $f=F(x^h)$. Both have order $n$, so $\xi^n=1$ and $h\mid n$.

The element $\alpha=\beta^h$ has order $e=n/h$. The polynomial $F$ is irreducible, since a proper factorization of $F$ would remain proper after substitution. Consequently $D=\ord_e(q)$ and
\[
 [\F_q(\beta):\F_q(\alpha)]=h.
\]
Thus $x^h-\alpha$ is irreducible over $\F_q(\alpha)$. If a prime $\ell$ divided $h$ but not $e$, choose an integer $u$ with $\ell u\equiv1\pmod e$. Then $(\alpha^u)^\ell=\alpha$, so $x^{h/\ell}-\alpha^u$ would be a proper factor of $x^h-\alpha$. This contradiction proves $\rad(h)\mid e$.

Now $\Phi_n(x)=\Phi_e(x^h)$. Each irreducible factor $G$ of $\Phi_e$ has degree $D$, so $G(x^h)$ divides $\Phi_n$ and has degree $hD=d$. All irreducible factors of $\Phi_n$ have degree $d$; hence $G(x^h)$ is irreducible. This gives the factor bijection. Substitution preserves weight, and the common support divisor in $\Phi_n$ is $h$, so every $G$ is a trinomial with support divisor one. The layer determines $h$ and hence $e=n/h$.
\end{proof}

\begin{corollary}[Lifting criterion]\label{cor:lift}
If $e>1$, $\gcd(e,q)=1$, $\rad(h)\mid e$, and
\begin{equation}\label{eq:lift}
 \ord_{eh}(q)=h\ord_e(q),
\end{equation}
then $\T_q(eh)$ is equivalent to $\T_q(e)$, and the irreducible factors correspond by $G(x)\mapsto G(x^h)$.
\end{corollary}
\begin{proof}
Use Lemma~\ref{lem:cyclosub} and compare factor degrees as in Proposition~\ref{prop:descent}.
\end{proof}

\begin{lemma}\label{lem:scalar}
Let $s\mid q-1$, $\gcd(s,m)=1$, and $\gcd(m,q)=1$. Put $d=\ord_m(q)$, and let $\R_s$ be the primitive $s$th roots in $\F_q$, with $\R_1=\{1\}$. Then
\begin{equation}\label{eq:scalar}
 \Phi_{sm}(x)=\prod_{\lambda\in\R_s}\ \prod_{G\mid\Phi_m\ \mathrm{monic\ irreducible}}
 \lambda^dG(x/\lambda).
\end{equation}
All displayed factors are distinct and irreducible. Scalar separation preserves weight and support divisor. For $m>1$, it preserves the property $\T_q$.
\end{lemma}
\begin{proof}
Every primitive $sm$th root has a unique decomposition $\lambda\alpha$, with $\lambda\in\R_s$ and $\alpha$ primitive of order $m$. Since Frobenius fixes $\lambda$, the minimal polynomial is $\lambda^dG(x/\lambda)$. The decomposition proves distinctness and completeness, and scalar substitution preserves the support. For $m=1$, this is the linear factorization of $\Phi_s$.
\end{proof}

\subsection{Local orders}
\begin{lemma}\label{lem:orders}
Let $\ell\nmid q$ be an odd prime, $d_\ell=\ord_\ell(q)$, and $b_\ell=v_\ell(q^{d_\ell}-1)$. Then
\begin{equation}\label{eq:odd-local}
 \ord_{\ell^a}(q)=d_\ell\ell^{\max\{a-b_\ell,0\}}\qquad(a\ge1).
\end{equation}
For odd $q$, put $\beta=v_2(q-1)$ and $\gamma=v_2(q+1)$. Then
\begin{equation}\label{eq:two-local}
 \ord_{2^a}(q)=\begin{cases}
 1,&a\le\beta,\\
 2,&\beta<a\le\beta+\gamma,\\
 2^{a-\beta-\gamma+1},&a>\beta+\gamma.
 \end{cases}
\end{equation}
Consequently, for any prime power $q$, if $\rad(n)\mid q^2-1$, and
$c=\gcd(n,q^2-1)$, $h=n/c$, then
\begin{equation}\label{eq:radial-lift}
 \ord_n(q)=h\ord_c(q),\qquad\ord_c(q)\in\{1,2\}.
\end{equation}
\end{lemma}
\begin{proof}
For odd $\ell$, the binomial theorem gives
\[
 v_\ell(q^{d_\ell t}-1)=b_\ell+v_\ell(t).
\]
Indeed, writing $q^{d_\ell}=1+\ell^{b_\ell}u$ with $\ell\nmid u$, raising to an exponent prime to $\ell$ leaves the valuation unchanged, while raising to $\ell$ increases it by one. Iteration proves the identity and then \eqref{eq:odd-local}.

For odd $q$ and odd $t$, $v_2(q^t-1)=\beta$. For even $t$, factoring the first square and then squaring successively gives
\[
 v_2(q^t-1)=\beta+\gamma+v_2(t)-1.
\]
Here every square of an odd integer is $1$ modulo $8$, so each subsequent factor of the form $u+1$ has valuation one. These two formulas imply \eqref{eq:two-local}.

For \eqref{eq:radial-lift}, every odd prime dividing $n$ has $d_\ell=1$ or $2$. The local excess beyond $v_\ell(q^2-1)$ contributes precisely the $\ell$-part of $h$ by \eqref{eq:odd-local}. If $q$ is odd, the $2$-part beyond $\beta+\gamma$ contributes the $2$-part of $h$ to a local order already containing a factor two. If $\ord_c(q)=1$, such a $2$-adic excess cannot occur, since it would force $2^{\beta+\gamma}\mid c$ but $c\mid q-1$. Taking least common multiples of the local orders now gives \eqref{eq:radial-lift}.
\end{proof}

\subsection{A degree bound from a missing prime}
\begin{lemma}[Missing-prime bound]\label{lem:primebound}
If $e\in\E_q$, $D=\ord_e(q)$, and $\ell$ is an odd prime not dividing $qe$, then $D\le\ell$.
\end{lemma}
\begin{proof}
Choose $F=x^D+A x^K+B$ with $\gcd(D,K)=1$, and a primitive $\ell$th root $\omega$. Set
\[
 R(U,V)=\prod_{j=0}^{\ell-1}(1+\omega^{jD}U+\omega^{jK}V).
\]
The coefficients of $1,U^\ell,V^\ell$ are all one. There is at least one more nonzero monomial: otherwise $R=1+U^\ell+V^\ell$, whereas the factor $1+U+V$ would force
\[
 1+U^\ell-(1+U)^\ell=0,
\]
whose coefficient of $U$ is $-\ell\ne0$.

If $\alpha$ is a root of $F$, then $\alpha^\ell$ still has order $e$. Its minimal polynomial $G$ is therefore another degree-$D$ trinomial in $\Phi_e$. Factoring over the algebraic closure gives
\[
 G(x^\ell)=\prod_{j=0}^{\ell-1}F(\omega^jx)
 =B^\ell R(x^D/B,A x^K/B).
\]
The leading scalar is one because $\ell$ is odd. If the exponents of distinct nonzero monomials of $R$ remained distinct after substitution, the right side would have at least four terms. Hence
\[
 (i-i')D+(j-j')K=0
\]
for two distinct monomials $U^iV^j,U^{i'}V^{j'}$. Necessarily $j\ne j'$, and coprimality gives $D\mid j-j'$. Since $R$ has total degree $\ell$, $0<|j-j'|\le\ell$, proving the bound.
\end{proof}

\begin{lemma}\label{lem:oddodd}
If $q$ and $e$ are odd and $e\in\E_q$, then $\ord_e(q)=2$.
\end{lemma}
\begin{proof}
Use the same argument with the square map. The two-variable product is
\[
 (1+U+V)(1+(-1)^D U+(-1)^K V).
\]
Since $D,K$ are not both even, it is one of
\[
 1-U^2-2UV-V^2,\qquad 1+2U+U^2-V^2,\qquad 1+2V+V^2-U^2.
\]
Each has four terms. Squaring preserves order $e$, so the collision argument gives $D\le2$. A trinomial has degree at least two.
\end{proof}

\section{Odd characteristic}\label{sec:odd}
For even root orders, squaring need not preserve the order. We instead use a quadratic coefficient norm, with a separate argument for the boundary of the $2$-primary part.

\begin{lemma}[Quadratic norm obstruction]\label{lem:norm2}
Let $q$ be odd, let $D>2$ be even, and let $0<K<D/2$ with $\gcd(D,K)=1$. If $C\in\F_q^*$ and $U^2\in\F_q^*$ is a nonsquare, then
\[
 W(x)=(x^D+C)^2-U^2x^{2K}
\]
cannot be the product of two degree-$D$ trinomials over $\F_q$ whose positive support exponents are coprime.
\end{lemma}
\begin{proof}
Suppose the factors are $x^D+u x^J+v$ and $x^D+w x^L+z$, with all four displayed lower coefficients nonzero. If $J\ne L$, the term of degree $D+\max(J,L)$ has a unique source and cannot vanish. This degree is larger than $D$, whereas $W$ has no term between degrees $D$ and $2D$. Thus $J=L$ and $w=-u$. The product is
\[
 x^{2D}+(v+z)x^D-u^2x^{2J}+u(z-v)x^J+vz.
\]
The case $2J>D$ gives an impossible term, and $2J=D$ contradicts $D>2$ and $\gcd(D,J)=1$. Hence $2J<D$. The largest nonzero term strictly below degree $D$ is then $-u^2x^{2J}$, so $J=K$. Comparison at degrees $J$ and $D$ yields $z=v$ and $2v=2C$. Finally, comparison at degree $2K$ gives $u^2=U^2$, a contradiction.
\end{proof}

\begin{theorem}\label{thm:odd-bases}
Over a finite field of odd characteristic, every base order has reduced degree two.
\end{theorem}
\begin{proof}
Let $e\in\E_q$, and write
\[
 e=2^a m,\qquad m\text{ odd},\qquad b=v_2(q-1),\qquad D=\ord_e(q).
\]
All irreducible factors of $\Phi_e$ have support divisor one by Lemma~\ref{lem:common-support}. If $a\le b$, scalar separation removes $2^a$. The case $m=1$ would give only linear factors and is impossible; otherwise the resulting odd-order layer is still a base layer, and Lemma~\ref{lem:oddodd} gives $D=2$.

Suppose $a>b$ and $D>2$. Then $D$ is even. Select a factor $f=x^D+A x^K+B$; replacing it by its reciprocal when necessary, assume $K<D/2$. The integer $K$ is odd.

\smallskip\noindent\emph{Case 1: $a\ge b+2$.}
Choose a primitive $2^{b+1}$st root $\eta$. Since $q\equiv1+2^b\pmod{2^{b+1}}$, we have $\eta^q=-\eta$. In particular, $\eta\in\F_{q^2}$ and $\eta^2\in\F_q$ is a nonsquare: its order has the full $2$-part $2^b$ of $q-1$.

Multiplication by $\eta$ or $-\eta$ preserves the order of every primitive $e$th root. Indeed, its $2$-primary exponent is odd, and the exponent contributed by a root of order at most $2^{b+1}$ is even when expressed modulo $2^a$. The odd-primary components are unchanged.

Set $C=B\eta^D$ and $U=A\eta^{D-K}$. Then $C\in\F_q^*$, $U^q=-U$, and $U^2$ is a nonsquare in $\F_q$: $D-K$ is odd. The two coefficient-conjugate polynomials
\[
 \eta^D f(x/\eta)=x^D+Ux^K+C,
 \qquad x^D-Ux^K+C
\]
have disjoint root sets, because their difference is $2Ux^K$ and their constant coefficient is nonzero. Each is separable. Their product $W=(x^D+C)^2-U^2x^{2K}$ is therefore squarefree, lies in $\F_q[x]$, and all of its roots are primitive $e$th roots. Every irreducible factor of $W$ has degree $D$, so there are exactly two. By $\T_q(e)$ and the common-support lemma they satisfy the hypotheses excluded by Lemma~\ref{lem:norm2}.

\smallskip\noindent\emph{Case 2: $a=b+1$.}
Fix a primitive $2^a$th root $\eta$, so again $\eta^q=-\eta$. Put $d_m=\ord_m(q)$. If $d_m$ is odd, then $D=2d_m$. For every primitive $m$th root $\beta$ and every primitive $2^a$th root $\eta'$, one has
\[
 (\eta'\beta)^{q^{d_m}}=-\eta'\beta.
\]
Thus every irreducible factor of $\Phi_e$ has a root set stable under negation. A monic polynomial of even degree with this property is an even polynomial, contradicting support divisor one. This argument also covers $m=1$ under the convention $d_m=1$.

We may therefore assume $d_m$ is even; then $D=d_m=2t$. Let $\beta$ be an arbitrary primitive $m$th root, and let $f=x^{2t}+A x^K+B$ be the minimal polynomial of $\eta\beta$. Here $K$ may depend on $\beta$, but $\gcd(2t,K)=1$, so $K$ is odd and $K\ne t$. Multiplying $f(x)$ and $f(-x)$, substituting $x^2=\eta^2y$, and dividing by $\eta^{2D}$ gives
\begin{equation}\label{eq:critical-shape}
 G(y)=(y^t+u)^2+v y^K,
 \qquad u=B\eta^{-D},\quad v=-A^2\eta^{2K-2D}.
\end{equation}
Both $u$ and $v$ belong to $\F_q^*$. This monic polynomial has degree $D$ and vanishes at $\beta^2$. Since $m$ is odd, $\beta^2$ still has order $m$ and degree $D$. Therefore $G$ is its minimal polynomial.

As $\beta$ varies, $\beta^2$ runs through every primitive $m$th root. Consequently every irreducible factor of $\Phi_m$ has the four-term form \eqref{eq:critical-shape}, and in particular satisfies
\begin{equation}\label{eq:middle-square}
 \bigl([y^t]G\bigr)^2=4G(0).
\end{equation}
Square the roots once more, and let $H$ be the resulting minimal polynomial. Its degree is again $D$, and $G(y)G(-y)=H(y^2)$.

If $t$ is even, then
\[
 H(z)=(z^{t/2}+u)^4-v^2z^K.
\]
Since $K\ne t$, its coefficient of $z^t$ is $6u^2$ and its constant coefficient is $u^4$. Applying \eqref{eq:middle-square} to $H$ gives $36u^4=4u^4$, hence $32u^4=0$, impossible in odd characteristic.

If $t$ is odd, direct multiplication gives
\[
 H(z)=z^{2t}-2u^2z^t-4uvz^{(t+K)/2}-v^2z^K+u^4.
\]
The five exponents are distinct: an equality involving $(t+K)/2$ with $t$ or $K$ would imply $K=t$, and $0<t,K<D$. All coefficients are nonzero. This contradicts the four-term form of every irreducible factor in $\Phi_m$. Both cases rule out $D>2$.
\end{proof}

\begin{proposition}[Odd-characteristic base orders]\label{prop:odd-base-list}
For odd $q$, the base orders are precisely $e\mid q^2-1$ with $e\nmid2(q-1)$.
\end{proposition}
\begin{proof}
By Theorem~\ref{thm:odd-bases}, it remains to consider degree two. If $\ord_e(q)=2$, a primitive $e$th root $\alpha$ has minimal polynomial
\[
 x^2-(\alpha+\alpha^q)x+\alpha^{q+1}.
\]
Its middle coefficient vanishes exactly when $\alpha^{q-1}=-1$. Under $e\nmid q-1$, this is equivalent to $e\mid2(q-1)$: the condition says that $\alpha^{q-1}$ has order two. Hence the condition for a trinomial is exactly the one stated. It applies simultaneously to all primitive $e$th roots.
\end{proof}

\section{Characteristic two: cubic base orders}\label{sec:cubic}
We first treat cubic layers whose root order is prime to three, then remove this restriction by a coefficient-norm argument.

\begin{lemma}\label{lem:cubic-cube}
Let $q=2^r$, suppose $3\nmid e$ and $\T_q(e)$, and let $x^3+Ax+B$ be an irreducible factor of $\Phi_e$. Then $A^3=B^2$. There exists $c\in\F_q^*$ such that $A=c^2$, $B=c^3$, and
\[
 e=7\ord(c),\qquad 3\nmid r.
\]
\end{lemma}
\begin{proof}
For a root $\alpha$, the polynomial with roots $\alpha^3,\alpha^{3q},\alpha^{3q^2}$ is
\begin{equation}\label{eq:cube-poly}
 y^3+By^2+(A^3+B^2)y+B^3.
\end{equation}
The elementary symmetric functions of the three cubes are $B$, $A^3+B^2$, and $B^3$. Since $3\nmid e$, the roots of \eqref{eq:cube-poly} have order $e$ and degree three, so it is a minimal polynomial in the same layer. Its leading, quadratic, and constant coefficients are nonzero; hence $A^3=B^2$.

With $c=B/A$, this gives $A=c^2$ and $B=c^3$. Thus the original polynomial is $c^3((x/c)^3+x/c+1)$. The binary cubic $x^3+x+1$ is irreducible, and its roots have order seven. It remains irreducible over $\F_{2^r}$ exactly when $3\nmid r$. In that case $7\nmid q-1$, so the orders of $c$ and of a primitive seventh root are coprime. The root order is therefore $7\ord(c)$.
\end{proof}

\begin{lemma}[The excess three-part in degree three]\label{lem:cubic-excess}
Let $q=2^r\equiv1\pmod3$ and $b=v_3(q-1)$. If $\ord_e(q)=3$ and $v_3(e)=b+1$, then $\T_q(e)$ fails.
\end{lemma}
\begin{proof}
Suppose otherwise. By taking a reciprocal, choose a factor $x^3+Ax+B$ whose trace coefficient is zero, and let $\alpha$ be a root. Let $\eta$ be the $3$-primary component of $\alpha$ in the cyclic group $\F_{q^3}^*$. Then $\eta$ has order $3^{b+1}$. The element $c=\eta^3$ belongs to $\F_q$ and is not a cube there, because its order has $3$-part $3^b$. Consequently $\F_q(\eta)=\F_{q^3}$, and $1,\eta,\eta^2$ is a basis.

Choose a primitive cube root $\omega\in\F_q$. The conjugates of $\eta$ are $\eta,\omega\eta,\omega^2\eta$. For
\[
 z=u_0+u_1\eta+u_2\eta^2\qquad(u_i\in\F_q),
\]
the trace and second elementary symmetric function of its three conjugates are
\begin{equation}\label{eq:cubic-symmetric}
 \Tr(z)=u_0,\qquad s_2(z)=u_0^2+c u_1u_2.
\end{equation}
These formulas follow by expansion using $1+\omega+\omega^2=0$ and $\eta^3=c$. Applying them to $\alpha$ gives $u_0=0$ and $A=c u_1u_2\ne0$.

Now set $\beta=\eta\alpha=c u_2+u_1\eta^2$. Multiplication by $\eta$ squares the $3$-primary component of $\alpha$ and leaves all other primary components unchanged. Thus $\beta$ still has order $e$ and degree three. Formula~\eqref{eq:cubic-symmetric} gives
\[
 \Tr(\beta)=c u_2\ne0,\qquad s_2(\beta)=c^2u_2^2\ne0.
\]
Its norm is nonzero as well. Its minimal polynomial has four nonzero terms, contradicting $\T_q(e)$.
\end{proof}

\begin{theorem}[Cubic base orders]\label{thm:cubic-base}
Formula~\eqref{eq:base3} holds. If $3\nmid r$ and $s\mid q-1$, then
\begin{equation}\label{eq:cubic-base-factors}
 \Phi_{7s}(x)=\prod_{\lambda\in\R_s}
 (x^3+\lambda^2x+\lambda^3)(x^3+\lambda x^2+\lambda^3),
\end{equation}
and all displayed factors are irreducible.
\end{theorem}
\begin{proof}
If $r$ is odd, $\ord_3(q)=2$, so a root order of degree three cannot be divisible by three. Lemma~\ref{lem:cubic-cube}, after reciprocation if needed, yields the assertion.

If $r$ is even, put $b=v_3(q-1)$. Since $e\mid q^3-1$, Lemma~\ref{lem:orders} gives $v_3(e)\le b+1$. Equality is excluded by Lemma~\ref{lem:cubic-excess}. All of the $3$-part of $e$ therefore lies in $q-1$ and can be removed by scalar separation. The remaining layer still has degree three and no factor three in its root order. Apply Lemma~\ref{lem:cubic-cube} and restore the scalar part. The result is $e=7s$ with $s\mid q-1$ and $3\nmid r$.

Conversely, the two binary factors of $\Phi_7$ remain irreducible when $3\nmid r$. Since $\gcd(7,s)=1$, scalar separation gives \eqref{eq:cubic-base-factors}; see also \cite[Lemma~2.8]{ChengII}.
\end{proof}

\subsection{Cubic coefficient norms}
\begin{lemma}[Cubic coefficient norm]\label{lem:norm3}
Let $q=2^r\equiv1\pmod3$, let $e\in\E_q$, and suppose $D=\ord_e(q)>3$. Then
\[
 v_3(e)\le v_3(q-1).
\]
\end{lemma}
\begin{proof}
Suppose $a=v_3(e)>b=v_3(q-1)$. Since $e\mid q^D-1$, the valuation formula in Lemma~\ref{lem:orders} gives $3\mid D$. Choose $f=x^D+A x^K+B$ with $\gcd(D,K)=1$; in particular $3\nmid K$.

Choose a primitive $3^{b+1}$st root $\eta\in\F_{q^3}$ and put $\gamma=\eta^3\in\F_q$. The element $\gamma$ is a noncube. If $a\ge b+2$, multiplication by $\eta$ preserves the order of every primitive $e$th root, because the change in its $3$-primary exponent is divisible by three. If $a=b+1$, choose $\eta$ to be the $3$-primary component of a selected root $\alpha$ of $f$. Relative to this generator, all conjugates of $\alpha$ have $3$-primary exponents congruent to $1$ modulo three, since $q\equiv1\pmod3$. Multiplication by $\eta$ changes that congruence to $2$, so their orders are again preserved. In either case multiplication by $\eta\omega^j$, $j=0,1,2$, also preserves these orders, where $\omega$ is a primitive cube root in $\F_q$: its additional exponent is a multiple of $3^b$.

Set
\[
 C=B\eta^D\in\F_q^*,\qquad U=A\eta^{D-K},\qquad
 U^3=A^3\gamma^{D-K}.
\]
Since $3\nmid D-K$, the element $U^3$ is a noncube in $\F_q$. The coefficient conjugates of $\eta^D f(x/\eta)$ have the form
\[
 x^D+C+U\omega^j x^K\qquad(j=0,1,2),
\]
up to a permutation of $j$. Their product is the coefficient norm
\begin{equation}\label{eq:norm3}
 W(x)=(x^D+C)^3+U^3x^{3K}\in\F_q[x].
\end{equation}
They are separable and pairwise coprime: subtracting two of them gives a nonzero scalar times $x^K$, and their constants are nonzero. All roots of $W$ have order $e$. Thus $W$ is the product of exactly three distinct degree-$D$ irreducible factors of $\Phi_e$.

All three factors are trinomials with support divisor one. Since $3\mid D$, the polynomial $W$ belongs to $\F_q[x^3]$. The substitution $x\mapsto\omega x$ therefore permutes its three monic irreducible factors. None can be fixed: a fixed factor $G=x^D+u x^J+v$ would require $3\mid J$, contrary to $\gcd(D,J)=1$. Hence they form one orbit of length three. Multiplying that orbit gives
\begin{equation}\label{eq:norm3-orbit}
 W(x)=(x^D+v)^3+u^3x^{3J}.
\end{equation}
Neither $3K$ nor $3J$ equals $D$ or $2D$. Indeed, any such equality, together with the appropriate support coprimality, would force $D=3$. Comparing the coefficients of $x^{2D}$ in \eqref{eq:norm3} and \eqref{eq:norm3-orbit} gives $v=C$. Subtracting the common cube then gives $J=K$ and $u^3=U^3$, contradicting the noncube property.
\end{proof}

\begin{corollary}\label{cor:even-r}
If $r$ is even, every base order over $\F_{2^r}$ has degree two or three.
\end{corollary}
\begin{proof}
Suppose $e\in\E_q$ and $D=\ord_e(q)>3$. Lemma~\ref{lem:primebound} forces $3\mid e$. Lemma~\ref{lem:norm3} shows that its entire $3$-part divides $q-1$. Remove that part by scalar separation. The resulting layer still has degree $D$ and support divisor one, but its root order is prime to three. Lemma~\ref{lem:primebound} then gives $D\le3$, a contradiction.
\end{proof}

\section{Characteristic two: quartic base orders}\label{sec:quartic}
A quartic trinomial in characteristic two has middle exponent one or three if it is irreducible: the case $x^4+Ax^2+B$ is a square. The two remaining supports are interchanged by reciprocation.

\begin{lemma}\label{lem:power7}
In characteristic two, if $\alpha^4+A\alpha+B=0$, then $z=\alpha^7$ satisfies
\begin{equation}\label{eq:power7}
 z^4+ABz^3+A^3(A^4+B^3)z+B^7=0.
\end{equation}
If $\alpha$ has order $m$, $\ord_m(q)=4$, and $7\nmid m$, the displayed polynomial is the minimal polynomial of $\alpha^7$ over $\F_q$.
\end{lemma}
\begin{proof}
Using $\alpha^4=A\alpha+B$ gives $z=B\alpha^3+A^2\alpha+AB$. Reduction in the basis $1,\alpha,\alpha^2,\alpha^3$ gives
\begin{align*}
 z^2&=AB^2\alpha^3+(A^4+B^3)\alpha^2+A^2B^2,\\
 z^3&=A^6\alpha^3+AB^4\alpha^2+(A^4B^2+B^5)\alpha,\\
 z^4&=A^3B^4\alpha^3+A^2B^5\alpha^2+(A^9+AB^6)\alpha
       +A^8B+A^4B^4+B^7.
\end{align*}
Substitution makes each of the four coefficients in \eqref{eq:power7} zero. If $7\nmid m$, the element $\alpha^7$ has order $m$ and degree four, so its monic quartic annihilator is its minimal polynomial.
\end{proof}

\begin{theorem}[Quartic base orders]\label{thm:quartic-base}
Formula~\eqref{eq:base4} holds. If $r$ is odd and $s\mid q-1$, then
\begin{equation}\label{eq:quartic-base-factors}
 \Phi_{15s}(x)=\prod_{\lambda\in\R_s}
 (x^4+\lambda^3x+\lambda^4)(x^4+\lambda x^3+\lambda^4),
\end{equation}
and all displayed factors are distinct and irreducible.
\end{theorem}
\begin{proof}
Suppose $\ord_e(q)=4$ and $\T_q(e)$. Since $q$ is even, the three integers $q-1,q+1,q^2+1$ are pairwise coprime. The divisibility $e\mid q^4-1$ therefore implies that
\[
 s=\gcd(e,q-1),\qquad m=e/s
\]
satisfy $\gcd(s,m)=\gcd(m,q-1)=1$ and $\ord_m(q)=4$. Scalar separation gives $\T_q(m)$. Furthermore $7\nmid m$: otherwise $\ord_7(q)$ would divide both $4$ and $\ord_7(2)=3$, hence would be one, contrary to $\gcd(m,q-1)=1$.

After taking a reciprocal, select a factor $x^4+Ax+B$ of $\Phi_m$. Equation~\eqref{eq:power7} gives another irreducible quartic in the same layer. Its leading, cubic, and constant coefficients are nonzero. The trinomial condition consequently forces $A^4=B^3$. With $c=B/A$, one has $c^3=A$, $c^4=B$, and
\[
 x^4+Ax+B=c^4\bigl((x/c)^4+x/c+1\bigr).
\]
The binary polynomial $x^4+x+1$ is irreducible: it has neither a root in $\F_2$ nor the quadratic factor $x^2+x+1$. Its roots have order $15$, since their degree excludes orders $1$ and $3$, and their minimal polynomial is not $\Phi_5$. It remains irreducible over $\F_{2^r}$ only when $r$ is odd. In that case $\gcd(15,q-1)=1$. A root $c\theta$ therefore has order $15\ord(c)$; since $\gcd(m,q-1)=1$, this implies $c=1$ and $m=15$. Hence $e=15s$.

Conversely, for odd $r$ the two binary quartics in $\Phi_{15}$ remain irreducible over $\F_q$. Scalar separation, using $\gcd(15,s)=1$, yields \eqref{eq:quartic-base-factors}.
\end{proof}

\begin{corollary}\label{cor:quartic-coeff}
Let $x^4+Ax+B\in\F_{2^r}[x]$ be irreducible and let $n$ be its root order. Then
\[
 \T_{2^r}(n)\quad\Longleftrightarrow\quad A^4=B^3.
\]
For an irreducible $x^4+Ax^3+B$, the condition is $B=A^4$.
\end{corollary}
\begin{proof}
The explicit factors in Theorem~\ref{thm:quartic-base} give necessity. Conversely, $A^4=B^3$ gives the scalar $c=B/A$ and reduces the polynomial to $x^4+x+1$, so its root order is $15\ord(c)$ and Theorem~\ref{thm:quartic-base} applies. For the other support, the monic reciprocal is $x^4+(A/B)x+1/B$; the condition becomes $A^4=B$.
\end{proof}

\section{Higher degrees in characteristic two}\label{sec:high}
Corollary~\ref{cor:even-r} settles even $r$. For odd $r$, we use the following consequence of Vishne's parity theorem.

\begin{lemma}[Vishne]\label{lem:swan8}
Let $q$ be even, $0<K<D$, and $A,B\in\F_q^*$. If $8\mid D$ and $\gcd(D,K)=1$, then $x^D+Ax^K+B$ is reducible over $\F_q$.
\end{lemma}
By \cite[Theorem~3.4]{Vishne}, the number of irreducible factors is even when $8\mid D$. Here $K$ is odd, so the polynomial is separable, and $D=2K$ is excluded by $\gcd(D,K)=1$ and $D>2$.

\begin{lemma}\label{lem:power17}
In characteristic two, suppose $\alpha^{12}+A\alpha^K+1=0$. If $z=\alpha^{17}$, then the following polynomial vanishes at $z$:
\begin{align}
 K=1:\quad H_1(z)={}&z^{12}+A^2z^{10}+A^6z^6+Az^5+A^3z^3+A^{10}z^2
                   +(A^{17}+A^5)z+1;\label{eq:H1}\\
 K=5:\quad H_5(z)={}&z^{12}+A^8z^8+A^{17}z^5+A^4z^4+A^2z^2+Az+1.
 \label{eq:H5}
\end{align}
If $A\ne0$, each polynomial has at least seven nonzero terms.
\end{lemma}
\begin{proof}
Let $R_j$ be the remainder of $x^j$ modulo $x^{12}+Ax^K+1$ in $\F_2[A][x]$. Starting with $R_j=x^j$ for $0\le j<12$, the recurrence
\begin{equation}\label{eq:identity-rec}
 R_{j+12}=A R_{j+K}+R_j\qquad(j\ge0)
\end{equation}
gives $H_K(x^{17})\equiv0\pmod{x^{12}+Ax^K+1}$ for $K=1,5$; the coefficient calculation is given in Appendix~\ref{app:identities}. If $A\ne0$, only the coefficient $A^{17}+A^5$ in $H_1$ can vanish. Its other seven coefficients, and all seven coefficients of $H_5$, are nonzero.
\end{proof}

\begin{theorem}\label{thm:odd-r-high}
If $r$ is odd, every base order over $\F_{2^r}$ has degree two, three, or four.
\end{theorem}
\begin{proof}
Suppose $e\in\E_q$ has degree $D>4$. Lemma~\ref{lem:primebound} forces $3\mid e$. Since $r$ is odd, $\ord_3(q)=2$, hence $D$ is even and $D\ge6$. The same lemma now forces $5\mid e$. Since $\ord_5(q)=4$, we obtain $4\mid D$. Lemma~\ref{lem:swan8} excludes $8\mid D$.

If $D>17$, the missing-prime bound also forces $17\mid e$. But $\ord_{17}(2)=8$ and $r$ is odd, so $\ord_{17}(q)=8$, which would imply $8\mid D$. Thus $D\le17$. The simultaneous conditions $D>4$, $4\mid D$, and $8\nmid D$ leave only $D=12$.

Since $r$ is odd, $q-1$ is odd and $3\nmid q-1$, so $\gcd(12,q-1)=1$. For every prime $\ell\mid q-1$, Lemma~\ref{lem:orders} therefore gives
\[
 v_\ell(q^{12}-1)=v_\ell(q-1)+v_\ell(12)=v_\ell(q-1).
\]
Consequently $s=\gcd(e,q-1)$ contains the full primary part of $e$ for every prime dividing $q-1$. Put $m=e/s$. Then $\gcd(s,m)=\gcd(m,q-1)=1$. Scalar separation preserves the base property and gives $\ord_m(q)=12$.

The norm of a primitive $m$th root belongs to both $\F_q^*$ and $\mu_m$. It is therefore one. Since the characteristic is two, every monic irreducible factor of $\Phi_m$ has constant coefficient one. Select
\[
 F=x^{12}+Ax^K+1.
\]
Its support divisor is one, so $K\in\{1,5,7,11\}$. Taking a reciprocal reduces to $K=1$ or $5$.

Now $17\nmid m$, since $\ord_{17}(q)=8$ does not divide $12$. If $\alpha$ is a root of $F$, then $\alpha^{17}$ is again primitive of order $m$ and has degree twelve. Lemma~\ref{lem:power17} gives its monic degree-twelve annihilator, which must consequently be its minimal polynomial. That polynomial has at least seven terms, contradicting $\T_q(m)$.
\end{proof}

\begin{proof}[Proof of Theorem~\ref{thm:bases}]
The odd-characteristic assertion is Proposition~\ref{prop:odd-base-list}. In characteristic two, every irreducible quadratic with nonzero constant coefficient is a trinomial, because a missing middle term would make it a square. This gives \eqref{eq:base2}. The cubic and quartic assertions are Theorems~\ref{thm:cubic-base} and~\ref{thm:quartic-base}. Corollary~\ref{cor:even-r} and Theorem~\ref{thm:odd-r-high} exclude all higher degrees.
\end{proof}

\section{Lifting and explicit factorizations}\label{sec:assembly}
We apply Corollary~\ref{cor:lift} to the base orders in Theorem~\ref{thm:bases}.

\begin{proposition}[Quadratic sources]\label{prop:quadratic-lifts}
For odd $q$, the layers arising from quadratic base orders are exactly \eqref{eq:odd-main}. For even $q$, they are exactly \eqref{eq:class2}. In both cases let
\[
 c=\gcd(n,q^2-1),\qquad h=n/c.
\]
Then $\ord_n(q)=2h$, and the complete irreducible factorization is
\begin{equation}\label{eq:quadratic-factors}
 \Phi_n(x)=\prod_{\alpha}
 \left(x^{2h}-(\alpha+\alpha^q)x^h+\alpha^{q+1}\right),
\end{equation}
where $\alpha$ runs through representatives of the $q$-Frobenius orbits of primitive $c$th roots. There are $\varphi(c)/2$ factors.
\end{proposition}
\begin{proof}
Suppose a layer comes from a quadratic base $e$, so $n=eh_0$, $\rad(h_0)\mid e$, and $\ord_n(q)=2h_0$. Then $\rad(n)\mid q^2-1$ and $e\mid c$. For odd $q$, $e\nmid2(q-1)$ implies $c\nmid2(q-1)$. For even $q$, the integers $q-1$ and $q+1$ are coprime. Since $e\mid q^2-1$ but $e\nmid q-1$, some prime of $e$ divides $q+1$ and not $q-1$. Hence $\rad(n)\nmid q-1$.

Conversely, the stated conditions make $c$ a quadratic base order by Theorem~\ref{thm:bases}. Lemma~\ref{lem:orders} gives $\ord_n(q)=2n/c=2h$. Moreover, $c$ and $n$ have the same prime divisors, so $\rad(h)\mid c$. Apply Corollary~\ref{cor:lift} to the quadratic minimal polynomials of the primitive $c$th roots. Their middle coefficients are nonzero by the base-order conditions. This proves \eqref{eq:quadratic-factors} and the factor count.
\end{proof}

\begin{proposition}[Cubic sources]\label{prop:cubic-lifts}
Let $q=2^r$. The layers arising from cubic base orders are precisely \eqref{eq:class3}. With
\[
 s=\gcd(M,q-1),\qquad L=7^{a-1}M/s,
\]
one has $\ord_n(q)=3L$, and
\begin{equation}\label{eq:cubic-lifts}
 \Phi_n(x)=\prod_{\lambda\in\R_s}
 (x^{3L}+\lambda^2x^L+\lambda^3)
 (x^{3L}+\lambda x^{2L}+\lambda^3).
\end{equation}
There are $2\varphi(s)$ irreducible factors.
\end{proposition}
\begin{proof}
A lift from a cubic base has $3\nmid r$ and $n=7^aM$, where $a\ge1$, $\gcd(M,7)=1$, and $\rad(M)\mid q-1$. Write the base as $7s_0$, with $s_0\mid q-1$. Since it divides $n$, $s_0\mid s=\gcd(M,q-1)$. Put $H=M/s$ and $v=v_7(r)$. Using $\ord_{7^a}(2)=3\cdot7^{a-1}$ and the identity $\ord_N(u^r)=\ord_N(u)/\gcd(\ord_N(u),r)$ gives
\[
 \ord_{7^a}(q)=3\cdot7^{\max\{a-1-v,0\}},\qquad\ord_M(q)=H.
\]
The first formula is recorded in \cite{ChengCorr}. Since $7\nmid q-1$, we have $\gcd(7,H)=1$, so taking the least common multiple gives
\[
 \ord_n(q)=3\cdot7^{\max\{a-1-v,0\}}\frac{H}{\gcd(3,H)}.
\]
For irreducible lifting from $7s_0$, this must equal $3\cdot7^{a-1}M/s_0$. Comparing $7$-adic valuations first gives $a=1$ or $v=0$. The remaining equality is $s_0=s\gcd(3,H)$. As $s_0\mid s$, it is equivalent to $s_0=s$ and $3\nmid H$.

These are precisely the conditions in \eqref{eq:class3}. Under them the order is $3L$, and $n=(7s)L$ with $\rad(L)\mid7s$. Substituting $x^L$ in \eqref{eq:cubic-base-factors} gives \eqref{eq:cubic-lifts}; each factor has the full degree $3L$ and is irreducible. See also \cite[Lemma~2.9]{ChengII}.
\end{proof}

\begin{proposition}[Quartic sources]\label{prop:quartic-lifts}
Let $q=2^r$. The layers arising from quartic base orders are precisely \eqref{eq:class4}. Put
\[
 s=\gcd(M,q-1),\qquad h=3^{a-1}5^{b-1}M/s.
\]
Then $\ord_n(q)=4h$ and
\begin{equation}\label{eq:quartic-lifts}
 \Phi_n(x)=\prod_{\lambda\in\R_s}
 (x^{4h}+\lambda^3x^h+\lambda^4)
 (x^{4h}+\lambda x^{3h}+\lambda^4).
\end{equation}
There are $2\varphi(s)$ irreducible factors.
\end{proposition}
\begin{proof}
A lift from a quartic base has odd $r$ and base order $15s_0$, $s_0\mid q-1$. It follows that $n=3^a5^bM$ with $a,b\ge1$, $\gcd(M,15)=1$, and $\rad(M)\mid q-1$. Put $s=\gcd(M,q-1)$, $u=v_3(r)$, and $v=v_5(r)$. We have $s_0\mid s$.

For odd $r$, the local orders and valuations are
\[
 \ord_3(q)=2,\quad\ord_5(q)=4,\quad
 v_3(q^2-1)=1+u,\quad v_5(q^4-1)=1+v.
\]
Thus
\[
 \ord_{3^a}(q)=2\cdot3^{\max\{a-1-u,0\}},\quad
 \ord_{5^b}(q)=4\cdot5^{\max\{b-1-v,0\}},\quad
 \ord_M(q)=M/s.
\]
Since $\gcd(M,30)=1$, the least common multiple is
\begin{equation}\label{eq:quartic-order}
 \ord_n(q)=4\cdot3^{\max\{a-1-u,0\}}
                   5^{\max\{b-1-v,0\}}M/s.
\end{equation}
Equating this with the required degree $4n/(15s_0)$ gives
\[
 s_0=s,\qquad \max\{a-1-u,0\}=a-1,\qquad
 \max\{b-1-v,0\}=b-1.
\]
The last two equalities are equivalent to $a=1$ or $3\nmid r$, and $b=1$ or $5\nmid r$. Conversely, these conditions make \eqref{eq:quartic-order} equal to $4h$. Substitution into the quartic base factorization gives the result by Corollary~\ref{cor:lift}.
\end{proof}

\begin{proof}[Proof of Theorem~\ref{thm:main}]
Proposition~\ref{prop:descent} reduces every admissible layer to a unique base order. Theorem~\ref{thm:bases} lists all such bases. Propositions~\ref{prop:quadratic-lifts}--\ref{prop:quartic-lifts} give their necessary and sufficient lifting conditions and therefore prove the classification.

In characteristic two, the alternatives are disjoint. Every prime in \textup{(Q)} divides $q^2-1$, whereas \textup{(C)} contains seven and $3\nmid r$ implies $7\nmid q^2-1$. In \textup{(F)}, $r$ is odd and $5\mid n$, but $5\nmid q^2-1$. Finally, \textup{(C)} with odd $r$ cannot contain five, since $5\nmid q-1$; if $r$ is even, \textup{(F)} is unavailable. This proves the asserted disjointness.
\end{proof}

\section{Consequences and examples}\label{sec:consequences}
\subsection{Comparison with global sparsity}
Following \cite{OR,ChengI,ChengII}, call $x^n-1$ \emph{globally $3$-sparse} if all of its monic irreducible factors have weight at most three.

\begin{corollary}\label{cor:global}
Suppose $\T_q(n)$ holds. Then $x^n-1$ is not globally $3$-sparse if and only if $q=2^r$ and $n$ satisfies the quartic alternative \eqref{eq:class4}.
\end{corollary}
\begin{proof}
First, if $\rad(n)\mid q^2-1$, take any divisor $d\mid n$ and set $c_d=\gcd(d,q^2-1)$, $h_d=d/c_d$. Lemma~\ref{lem:orders} and cyclotomic substitution imply that every irreducible factor of $\Phi_d$ is obtained from a linear or quadratic factor of $\Phi_{c_d}$ by substituting $x^{h_d}$. It therefore has at most three terms. This covers all odd-characteristic cases and alternative \textup{(Q)}.

In alternative \textup{(C)}, write a divisor of $n$ as $7^iT$ with $T\mid M$. If $i=0$, then $\rad(T)\mid q-1$, and the same lifting argument in even characteristic gives binomial factors. If $i\ge1$, the conditions of \eqref{eq:class3} are inherited: $3\nmid T/\gcd(T,q-1)$ follows by comparing $3$-adic valuations, and $i=1$ or $7\nmid r$ follows from the corresponding condition on $a$. Proposition~\ref{prop:cubic-lifts} supplies trinomial factors. Hence every divisor layer has weight at most three.

In alternative \textup{(F)}, $r$ is odd and $5\mid n$. The equality $\ord_5(2^r)=4$ makes $\Phi_5=x^4+x^3+x^2+x+1$ irreducible over $\F_q$. It divides $x^n-1$ and has five terms, giving the opposite conclusion.
\end{proof}

\subsection{Small fields and lifting restrictions}
\begin{example}[Small fields]\label{ex:small}
The complete lists for the first three characteristic-two fields and for $\F_3$ are
\[
\begin{array}{c|l}
\F_q&\text{all admissible }n>1\\ \hline
\F_2&3^a,\ 7^a,\ 3^a5^b\quad(a,b\ge1),\\
\F_4&3^a5^b\ (a\ge0,b\ge1),\quad3^\varepsilon7^b\ (\varepsilon=0,1,b\ge1),\\
\F_8&3^a7^b\ (a\ge1,b\ge0),\quad3\cdot5^b7^c\ (b\ge1,c\ge0),\\
\F_3&2^a\quad(a\ge3).
\end{array}
\]
For $\F_3$, the radical condition allows only powers of two, and $\gcd(2^a,8)\nmid4$ is equivalent to $a\ge3$. Writing $h=2^{a-3}$ gives the complete factorization
\[
 \Phi_{2^a}(x)=(x^{2h}+x^h-1)(x^{2h}-x^h-1).
\]
The two base quadratics are irreducible over $\F_3$, and $\ord_{2^a}(3)=2h$ verifies every subsequent lift.

Over $\F_5$, the admissible orders are instead $n=2^a3^b$ with $a\ge0$ and $b\ge1$. For example,
\[
 \Phi_{12}(x)=(x^2+2x-1)(x^2-2x-1)
\]
is a trinomial irreducible factorization. The radical condition alone is insufficient: $\rad(8)\mid24$, but $\gcd(8,24)\mid2(5-1)$, and
$\Phi_8(x)=(x^2+2)(x^2-2)$ has two irreducible binomial factors.
\end{example}

\begin{example}[Failure of lifting]\label{ex:saturation}
The restrictions at $3$, $5$, and $7$ in \eqref{eq:class3} and \eqref{eq:class4} cannot be omitted. In each row below, $e$ is a base order and $n=eh$, but the degree condition \eqref{eq:lift} fails and $\T_q(n)$ does not hold:
\begin{center}
\begin{tabular}{rrrrcc}
\toprule
$q$ & $n$ & $e$ & $h$ & $\ord_n(q)$ & $h\ord_e(q)$\\
\midrule
$4$   & $63$ & $21$ & $3$ & $3$ & $9$\\
$8$   & $45$ & $15$ & $3$ & $4$ & $12$\\
$32$  & $75$ & $15$ & $5$ & $4$ & $20$\\
$128$ & $49$ & $7$  & $7$ & $3$ & $21$\\
\bottomrule
\end{tabular}
\end{center}
For the first row, let $\F_4=\F_2(\omega)$ with $\omega^2+\omega+1=0$. The polynomial $P=x^3+x^2+x+\omega$ has no root in $\F_4$, so it is irreducible. A root $\beta$ satisfies $\beta^{21}=\omega$, hence has order $9$ or $63$. Since $\Phi_9=(x^3+\omega)(x^3+\omega^2)$, its order cannot be nine. Thus $P$ is a four-term factor of $\Phi_{63}$; see \cite[Proposition~3.3]{ChengII}.

Over both $\F_2$ and $\F_{32}$,
\[
 \Phi_{45}(x)=(x^{12}+x^3+1)(x^{12}+x^9+1)
\]
is an irreducible factorization. Over $\F_8=\F_2(u)$, where $u^3+u+1=0$, the five-term polynomial
\[
 x^4+u^2x^3+(u+1)x^2+ux+1
\]
divides $x^{12}+x^3+1$. Its roots have order $45$, and its degree is $\ord_{45}(8)=4$, so it is irreducible over $\F_8$. The second and third rows reflect the separate conditions $3\nmid r$ and $5\nmid r$ for quartic lifting. An explicit four-term factor for the last row is given in \cite{ChengCorr}.
\end{example}

\subsection{Quartic factor counts}
For $\ord_n(q)=4$, let $N_j(n;q)$ denote the number of monic irreducible factors of $\Phi_n$ with exactly $j$ nonzero terms.

\begin{theorem}[Quartic gap]\label{thm:gap}
Let $q=2^r$ and $\ord_n(q)=4$. If $\T_q(n)$ fails, then
\begin{equation}\label{eq:gap}
 N_3(n;q)\le N_4(n;q),\qquad
 \frac{N_3(n;q)}{\varphi(n)/4}\le\frac12.
\end{equation}
Thus the trinomial proportion in a quartic layer is either one or at most one half.
\end{theorem}
\begin{proof}
As in the quartic classification, write $n=sm$, where $s=\gcd(n,q-1)$, $\gcd(s,m)=\gcd(m,q-1)=1$, and $7\nmid m$. Scalar separation gives
\[
 N_j(n;q)=\varphi(s)N_j(m;q)
\]
for every $j$. It suffices to work at order $m$.

The map taking each root to its seventh power permutes the irreducible factors of $\Phi_m$, because it is a bijection on the primitive $m$th roots and commutes with Frobenius. If $x^4+Ax+B$ is a trinomial factor, then $A^4\ne B^3$. Otherwise the argument in Theorem~\ref{thm:quartic-base} would force $m=15$ and make the entire layer trinomial. By Lemma~\ref{lem:power7}, its image is therefore a four-term factor. A factor with middle exponent three has the same property: reciprocation preserves weight and commutes with the seventh-power map. Hence this permutation sends all trinomial factors injectively into the set of four-term factors. This proves the first inequality. The second follows from $N_3+N_4\le\varphi(m)/4$, and then from scalar separation.
\end{proof}

\begin{proposition}\label{prop:quartic-count}
Let $q=2^r$. If $r$ is even, there is no irreducible quartic trinomial over $\F_q$. If $r$ is odd and $A,B\in\F_q^*$, let $c\in\F_q^*$ be the unique element with $c^3=A$. Then
\begin{equation}\label{eq:trace-quartic}
 x^4+Ax+B\text{ is irreducible}\quad\Longleftrightarrow\quad
 \Tr_{\F_q/\F_2}(B/c^4)=1.
\end{equation}
For odd $r$, there are $q(q-1)$ monic irreducible quartic trinomials, of which precisely $2(q-1)$ belong to layers satisfying $\T_q$. Their proportion is $2/q$.
\end{proposition}
\begin{proof}
For odd $r$, the cube map permutes $\F_q^*$, so scaling reduces to $y^4+y+t$. On $\F_q$, the map $y\mapsto y^4+y$ has kernel $\F_2$ and image the absolute-trace-zero subspace: its image has $q/2$ elements and is contained in that subspace. Thus trace zero gives a root in $\F_q$.

If $\Tr_{\F_q/\F_2}(t)=1$, there is no root in $\F_{q^2}$ either. Indeed, taking the trace of $y^4+y=t$ from $\F_{q^2}$ to $\F_4$ gives zero on the left, while on the right
\[
 \sum_{i=0}^{r-1}t^{4^i}=\Tr_{\F_q/\F_2}(t)=1,
\]
since multiplication by two permutes the indices modulo odd $r$. A reducible quartic must have a factor of degree one or two, so the polynomial is irreducible. Each $A\ne0$ admits $q/2$ choices of $t$, giving $q(q-1)/2$ trinomials with middle exponent one. Reciprocation gives the same number with exponent three; these sets are disjoint.

For even $r$, take a cube root $c$ of $A$ in $\F_{q^3}$ and scale again. This field contains $\F_4$. Put $Q=q^3$. For a root $\beta$ of $y^4+y+t$, the difference $\beta^Q-\beta$ belongs to $\F_4$, since it is in the kernel of $y^4+y$. Frobenius of order $Q$ fixes that difference, so $\beta^{Q^2}=\beta$. There can be no irreducible quartic over $\F_Q$. But an irreducible quartic over $\F_q$ would remain irreducible over a cubic extension, a contradiction. This nonexistence is also a special case of \cite[Corollary~5.1]{Vishne}.

Finally, Theorem~\ref{thm:quartic-base} gives $2\varphi(s)$ factors for each $s\mid q-1$. Distinct root orders cannot share an irreducible factor. Therefore their total is
\[
 \sum_{s\mid q-1}2\varphi(s)=2(q-1).
\]
Dividing by $q(q-1)$ proves the last assertion.
\end{proof}

\begin{remark}[Indices divisible by the characteristic]\label{rem:p-part}
For $n=p^a m$ with $a\ge1$ and $p\nmid m$, reduction modulo $p$ gives
\[
 \Phi_{p^a m}(x)=\Phi_m(x)^{p^{a-1}(p-1)}.
\]
This follows by applying Frobenius to $\Phi_{p^a m}(x)=\Phi_m(x^{p^a})/\Phi_m(x^{p^{a-1}})$. Thus the distinct factors are those of $\Phi_m$, and only their multiplicities change. For $m>1$, Theorem~\ref{thm:main} applies; for $m=1$, the only distinct factor is $x-1$.
\end{remark}

\appendix
\section{Polynomial identities used in the degree-twelve exclusion}\label{app:identities}
Let $\mathbf r_j\in\F_2[A]^{12}$ be the coefficient vector of $R_j$ from \eqref{eq:identity-rec}, with $\mathbf r_0,\ldots,\mathbf r_{11}$ the coordinate vectors. The recurrence
\[
 \mathbf r_{j+12}=A\mathbf r_{j+K}+\mathbf r_j
\]
reduces $H_K(x^{17})$ modulo $x^{12}+Ax^K+1$. Writing $H_K(z)=\sum_i c_i(A)z^i$, its remainder vector is $\sum_i c_i(A)\mathbf r_{17i}$. The coefficients are listed below; blank entries are zero.
\begin{center}
\begin{tabular}{c|cccccccccc}
\toprule
$i$&12&10&8&6&5&4&3&2&1&0\\
\midrule
$K=1$&$1$&$A^2$&&$A^6$&$A$&&$A^3$&$A^{10}$&$A^{17}+A^5$&$1$\\
$K=5$&$1$&&$A^8$&&$A^{17}$&$A^4$&&$A^2$&$A$&$1$\\
\bottomrule
\end{tabular}
\end{center}
The following Python code computes the two zero remainder vectors. An integer encodes a polynomial in $\F_2[A]$ by its binary digits, so exclusive-or is addition and a shift by $a$ places multiplies by $A^a$. The tuple at index $i$ lists the exponents of $A$ occurring in $c_i(A)$.
\par\medskip
\noindent\begin{minipage}{\linewidth}
\small
\begin{verbatim}
terms = {
    1: {12: (0,), 10: (2,), 6: (6,), 5: (1,),
        3: (3,), 2: (10,), 1: (17, 5), 0: (0,)},
    5: {12: (0,), 8: (8,), 5: (17,), 4: (4,),
        2: (2,), 1: (1,), 0: (0,)}
}
for K, H in terms.items():
    R = [[int(i == j) for i in range(12)] for j in range(12)]
    for j in range(17 * max(H) - 11):
        R.append([(a << 1) ^ b for a, b in zip(R[j+K], R[j])])
    rem = [0] * 12
    for i, powers in H.items():
        for a in powers:
            for j in range(12):
                rem[j] ^= R[17*i][j] << a
    assert rem == [0] * 12, K
\end{verbatim}
\end{minipage}
\par\medskip
Since $\deg H_1=\deg H_5=12$, the recurrence computes all remainders through $R_{204}$. Both assertions give
\[
 H_K(x^{17})\equiv0\pmod{x^{12}+Ax^K+1},\qquad K\in\{1,5\}.
\]


\begin{thebibliography}{99}
\bibitem{BGO}
F.~E. Brochero Mart\'{i}nez, C.~R. Giraldo Vergara, and L.~B. de Oliveira,
Explicit factorization of $x^n-1\in\F_q[x]$,
\emph{Designs, Codes and Cryptography} \textbf{77} (2015), 277--286.
\doi{10.1007/s10623-014-0005-y}.

\bibitem{ChengI}
K. Cheng, The $3$-sparsity of $X^n-1$ over finite fields,
\href{https://arxiv.org/abs/2507.06655v2}{arXiv:2507.06655v2}, 2025.

\bibitem{ChengII}
K. Cheng, The $3$-sparsity of $X^n-1$ over finite fields of characteristic $2$,
\emph{Transactions of the London Mathematical Society} \textbf{13} (2026), e70026.
\doi{10.1112/tlm3.70026}.
Corrected author version: \emph{The $3$-sparsity of $X^n-1$ over finite fields, II},
\href{https://arxiv.org/abs/2507.10779v2}{arXiv:2507.10779v2}.

\bibitem{ChengCorr}
K. Cheng, Corrigendum: The $3$-sparsity of $X^n-1$ over finite fields of characteristic $2$,
\emph{Transactions of the London Mathematical Society} \textbf{13} (2026), e70035.
\doi{10.1112/tlm3.70035}.

\bibitem{Graner}
A.-M. Graner, Closed formulas for the generators of all constacyclic codes and for the factorization of $X^n-1$, the $n$-th cyclotomic polynomial and every composition of the form $f(X^n)$ over a finite field for arbitrary positive integers $n$,
\emph{Finite Fields and Their Applications} \textbf{109} (2026), 102695.
\doi{10.1016/j.ffa.2025.102695}.
Preprint version: \href{https://arxiv.org/abs/2306.11183}{arXiv:2306.11183}.

\bibitem{LN}
R. Lidl and H. Niederreiter, \emph{Finite Fields}, 2nd ed., Encyclopedia of Mathematics and its Applications, vol.~20, Cambridge University Press, Cambridge, 1997.
\doi{10.1017/CBO9780511525926}.

\bibitem{OR}
D. Oliveira and L. Reis, On polynomials $x^n-1$ over binary fields whose irreducible factors are binomials and trinomials,
\emph{Finite Fields and Their Applications} \textbf{73} (2021), 101837.
\doi{10.1016/j.ffa.2021.101837}.

\bibitem{Swan}
R.~G. Swan, Factorization of polynomials over finite fields,
\emph{Pacific Journal of Mathematics} \textbf{12} (1962), 1099--1106.
\doi{10.2140/pjm.1962.12.1099}.

\bibitem{Vishne}
U. Vishne, Factorization of trinomials over Galois fields of characteristic $2$,
\emph{Finite Fields and Their Applications} \textbf{3} (1997), 370--377.
\doi{10.1006/ffta.1997.0191}.

\bibitem{WW}
L. Wang and Q. Wang, On explicit factors of cyclotomic polynomials over finite fields,
\emph{Designs, Codes and Cryptography} \textbf{63} (2012), 87--104.
\doi{10.1007/s10623-011-9537-6}.

\bibitem{Zhu} L. Zhu, J. Liu, H. Wu, Explicit representatives and sizes of cyclotomic cosets and their application to cyclic codes over finite fields, Finite Fields and Their Applications, 111(2026), 102761.
\end{thebibliography}
\end{document}